\documentclass[a4paper,11pt]{article}

\usepackage[utf8]{inputenc}
\usepackage[titletoc,toc,title]{appendix}
\usepackage{a4wide}
\usepackage{authblk}

\usepackage{graphicx}

\usepackage{tikz}

\usepackage{booktabs}
\usepackage{array}
\usepackage{paralist}
\usepackage{verbatim}
\usepackage{subfigure}
\usepackage{algorithm}
\usepackage{algpseudocode}

\usepackage{amsmath,amssymb,amsfonts}
\usepackage{eufrak} 
\usepackage{mathrsfs}  
\usepackage{bm}  
\renewcommand{\theta}{\vartheta}
\renewcommand{\phi}{\varphi}
\renewcommand{\rho}{\varrho}
\renewcommand{\epsilon}{\varepsilon}

\renewcommand{\mathbf}{\boldsymbol}
\makeatletter
\newcommand{\pushright}[1]{\ifmeasuring@#1\else\omit\hfill$\displaystyle#1$\fi\ignorespaces}
\newcommand{\pushleft}[1]{\ifmeasuring@#1\else\omit$\displaystyle#1$\hfill\fi\ignorespaces}
\makeatother
\DeclareMathOperator{\argmin}{arg\,min}

\DeclareMathOperator{\ddiv}{div}

\newcommand*{\Gr}{\mathrm{Gr}\,}
\newcommand*{\de}{\mathop{}\!\mathrm{d}}
\newcommand*{\Proj}{\mathit{\Pi}}
\newcommand*{\Om}{\mathit{\Omega}}

\newenvironment{sistema}%
{\left\lbrace\begin{array}{@{}l@{}}}%
{\end{array}\right.}
\usepackage{footnote}

\usepackage{fancyhdr}
\fancypagestyle{FirstPage}{
\lfoot{The final publication is available at \url{link.springer.com}}
}

\usepackage[natbib=true,backend=bibtex]{biblatex}
\bibliography{rbsissa}

\begin{document}

\title{On the application of reduced basis methods to bifurcation problems in incompressible fluid dynamics}
\author[1]{Giuseppe Pitton}
\author[1]{Gianluigi Rozza}
\affil[1]{SISSA, International School for Advanced Studies, Mathematics Area, mathLab, via Bonomea 265, 34136, Trieste, Italy\\ Email: \texttt{gpitton@sissa.it, grozza@sissa.it}}
\date{}

\maketitle

\thispagestyle{FirstPage}

\begin{abstract}
In this paper we apply a reduced basis framework for the computation of flow bifurcation (and stability) problems in fluid dynamics. The proposed method aims at reducing the complexity and the computational time required for the construction of bifurcation and stability diagrams. The method is quite general since it can in principle be specialized to a wide class of nonlinear problems, but in this work we focus on an application in incompressible fluid dynamics at low Reynolds numbers. The validation of the reduced order model with the full order computation for a benchmark cavity flow problem is promising.

\end{abstract}

\section{Introduction}
Among the many applications of bifurcation theory, there is the computation of solution branches for nonlinear problems~\cite{Prodi}, and their intersection points (called \emph{bifurcation points}) as a function of some physical or geometrical control parameters.
Bifurcation diagrams show a qualitative or quantitative (pointwise or integral) aspect of a solution as a function of the control parameters, allowing to understand the structure and the physical features of the solution set.
See for example~\cite{Chandra} for a classical introduction to stability problems in fluid mechanics.


A classical application of bifurcation theory is the analysis of stability problems in Elasticity, where it is well known that for many structures there exists a \emph{critical load} which if exceeded will produce a catastrophic collapse (\cite{Timoshenko},~\cite{Marsden}).

Common applications of bifurcation theory in fluid mechanics include a wide range of industrial problems such as tribology, microfluid dynamics, biomedical industry, biomedicine (blood flows), and many more.
In particular, the field of hydrodynamic stability focuses on the classification of the nature of flows 
as some control parameters are varied~\cite{Chandra}. Specifically, the purpose of a stability investigation is, given a well-defined initial or boundary value problem, to classify the asymptotic solutions as stable, or time periodic, or even chaotic in time.

The construction of stability maps and bifurcation diagrams is a delicate and very expensive task, requiring an important computational effort and prohibitive if the number and range of parameters is large. This is particularly true for three and higher dimensional simulations.

Recent developments of reduced basis methods have focused on the reduction of computational time for a wide range of differential problems, while mantaining a prescribed tolerance on error bounds~\cite{Rozza:ARCME}. It is therefore of great interest to further investigate how such methods can be applied to stability problems in fluid dynamics and elasticity in order to reduce the required computational power.

The literature has already shown the effectiveness of the Proper Orthogonal Decomposition (POD) both for analysis of principal modes~\cite{Holmes:2012} and for the reduction of computational power required by transient simulations. For instance, Terragni and Vega~\cite{Terragni:2012} showed how a POD approach could save a considerable amount of computing time for the analysis of bifurcations in some nonlinear dissipative systems.
In a recent paper, Herrero, Maday and Pla~\cite{Maday:RB2} have shown that both POD and the Reduced Basis Method (RBM) can reconstruct the behaviour of velocity and temperature field for a two-dimensional natural convection (Boussinesq) problem with large reduction of the computational power with respect to classical techniques. In particular, stable and unstable solutions are correctly identified, and a surrogate error estimate is always mantained below a prescribed tolerance.
A recent remarkable work by Yano et al.~\cite{Yano:2013} introduced a RB Method for the stability of flows under perturbations in the forcing term or in the boundary conditions, based on a space-time framework that allows for particularly sharp error estimates.

Summarizing, given the relatively fast decay of the energy spectrum for flows at sufficiently low Reynolds numbers, a reduced basis method is expected to be an efficient tool for flow stability analysis.
Reduced basis techniques~\cite{AQGR:2014} historically have been proven effective for the study of elastic stability of plates~\cite{Noor:1983}, but their application in more complex parametrized stability problems such as in Fluid Mechanics is still an open and quite relevant research field. 

In this work, we propose the application of reduced basis techniques to cut the quite demanding computational costs for bifurcation detection and stability analysis of viscous flows. 
We will consider a benchmark problem well known in the literature, namely a buoyancy-driven flow in a rectangular cavity~\cite{Roux:GAMM}. 


The structure of the manuscript is the following. In Section~\ref{sec:appfluid} we introduce a typical model problem for bifurcating flows that will be used to fix the ideas in the exposition and later as a numerical test, together with its Galerkin approximation. In Section~\ref{sec:ROM} the reduced-order method used in this work for the approximation of the Rayleigh-B\'{e}nard equations is presented. In Section~\ref{sec:bifurcation_detection} we discuss a RB technique adopted for the detection of bifurcation points of laminar flows. The results of some numerical tests are shown and discussed in Section~\ref{sec:results}.

\section{Problem statement}
\label{sec:appfluid}

In this section we present a typical model problem for bifurcating flows in incompressible fluid dynamics. The model problem consists of a Rayleigh-B\'{e}nard cavity flow, introduced in~\cite{Roux:GAMM} and widely studied over the 90s. 
The particular instance of Rayleigh-B\'{e}nard cavity flow chosen as test case will prove to be a good candidate for our analysis, due to the richness of bifurcating behaviour despite the very simple geometrical setting.



\subsection{Mathematical model}

The strong formulation of the Rayleigh-B\'{e}nard cavity problem, as stated in~\cite{Roux:GAMM}, is given in adimensional variables as follows: find a velocity field $\mathbf{u}$ and a pressure field $p$ such that
\begin{equation}
  \begin{cases}
    \frac{\partial\mathbf{u}}{\partial t}+\mathbf{u}\cdot\nabla\mathbf{u}-\Delta\mathbf{u}+\nabla p=\Gr x\mathbf{\jmath} \qquad&\text{on }\mathit{\Omega}(A) \\
    \ddiv\mathbf{u}=0 &\text{on }\mathit{\Omega}(A) \\
    \mathbf{u}=0 & \text{on }\mathit{\Gamma}_D(A), \\
  \end{cases}
  \label{eq:strongRB}
\end{equation}
where $\mathbf{\jmath}$ is the unit vector directed along the vertical axis, $x$ the horizontal coordinate, $x$ the horizontal coordinate, $\Gr$  the \emph{Grashof number} wich expresses roughly the ratio of buoyancy to viscous forces, and is a parameter. The domain $\mathit{\Omega}(A)$ considered in the benchmark is a rectangular bidimensional cavity with unit height and width $A$, and we choose the bottom left corner as the origin of the coordinate system. We consider the problem with fully Dirichlet boundary conditions, in symbols we write $\mathit{\Gamma}_D(A)\equiv\partial\mathit{\Omega}(A)$.

Some representative solutions to equation~\eqref{eq:strongRB} for some values of the parameters $A$ and $\Gr$ are shown in Figure~\ref{fig:snapsGAMM}.
The flow consists of a small number (1 to 5) of vortices, and for some values of the parameters it may happen that more than one flow configurations are possible.
In particular, for small values of $\Gr$ only the solution with one vortex appears, regardless of the cavity length $A$. However, if for a fixed value of $A$ the Grashof number is slowly increased, then at some point a flow pattern with more than one vortex will appear (the precise number of vortices depends on $A$).



The well-posedness of system~\eqref{eq:strongRB} can be proved (see e.g.~\cite{QuarteroniValli}) by considering its weak formulation. Introducing the velocity and pressure spaces $\mathbf{V}\equiv[H^1_0(\mathit{\Omega(A)})]^d$ and $Q\equiv L^2_0(\mathit{\Omega(A)})$, respectively\footnote{The notation chosen for the function spaces may need clarification. We refer with $H^1_0$ to the Sobolev space with zero trace at the boundary (velocity), and with $L^2_0$ to the Lebesgue $L^2$ functions with zero average (pressure).}, where $d$ is the spatial dimension. Then, multiplying the momentum and mass balance equations in~\eqref{eq:strongRB} respectively by the test functions $\mathbf{v}\in\mathbf{V}$ and $q\in Q$, and integrating formally by parts, we get the variational formulation that reads as follows:
find $(\mathbf{u},p)\in \mathbf{V}\times Q$ such that
\begin{equation}
  \begin{cases}
    m(\mathbf{u},\mathbf{v})+c(\mathbf{u},\mathbf{u},\mathbf{v})+a(\mathbf{u},\mathbf{v})+b(\mathbf{v},p)=f(\mathbf{v})\qquad &\forall\mathbf{v}\in \mathbf{V} \\
    b(\mathbf{u},q)=0 &\forall q\in Q,
  \end{cases}
  \label{eq:weakNS}
\end{equation}
where the following bilinear forms have been introduced
\begin{equation}
  \begin{aligned}
    &a(\mathbf{u},\mathbf{v})=\int_{\mathit{\Omega(A)}}\nabla\mathbf{v}:\nabla\mathbf{u}\de\mathbf{x}\\
    &b(\mathbf{v},q)=\int_{\mathit{\Omega(A)}}q\ddiv\mathbf{v}\de\mathbf{x}\\
    &m(\mathbf{u},\mathbf{v})=\int_{\mathit{\Omega(A)}}\mathbf{v}\cdot\frac{\partial\mathbf{u}}{\partial t}\de\mathbf{x},
  \end{aligned}
  \label{eq:bilinearforms}
\end{equation}
along with the variational form
\begin{equation}
  c(\mathbf{u},\mathbf{w},\mathbf{v})=\int_{\mathit{\Omega(A)}}\mathbf{v}\cdot\left(\mathbf{u}\cdot\nabla\mathbf{w}\right)\de \mathbf{x}
  \label{eq:formC}
\end{equation}
and the linear form
\begin{equation}
        f(\mathbf{v})=\int_{\mathit{\Omega(A)}}\Gr\theta\mathbf{\jmath}\cdot\mathbf{v}\de\mathbf{x}.
  \label{eq:defF}
\end{equation}

\subsection{Parametrized formulation}
We are interested in approximating the solutions set of equations~\eqref{eq:weakNS} for a wide range of aspect ratios for the rectangular geometry, and for an interval of the Grashof number over which the bifurcations take place. Specifically, the independent parameters for this problem are the cavity length $\mu\equiv A$ and the Grashof number $\Gr$. To simplify the notation, let us introduce the parameter vector $\mathbf{\mu}=(\mu,\Gr)$.

In the following, the dependence of the variational forms on the parameter vector $\mathbf{\mu}$ will be denoted explicitly, and problem~\eqref{eq:weakNS} will be cast as:
Find $(\mathbf{u}(\mathbf{\mu}),p(\mathbf{\mu}))\in \mathbf{V}\times Q$ such that
\begin{equation}
  \begin{split}
  \begin{cases}
    m(\mathbf{\mu};\mathbf{u}(\mathbf{\mu}),\mathbf{v})+c(\mathbf{\mu};\mathbf{u}(\mathbf{\mu}),\mathbf{u}(\mathbf{\mu}),\mathbf{v})+a(\mathbf{\mu};\mathbf{u}(\mathbf{\mu}),\mathbf{v})+b(\mathbf{\mu};\mathbf{v},p(\mathbf{\mu}))=&f(\mathbf{\mu};\mathbf{v})  \\ &\forall\mathbf{v}\in \mathbf{V} \\
    b(\mathbf{\mu};\mathbf{u}(\mathbf{\mu}),q)=0 & \forall q\in Q.
  \end{cases}
  \end{split}
  \label{eq:weakRBmu}
\end{equation}
We remark that each of the operators can be written as the product of a part which depends on $\mathbf{\mu}$ and a part depending on $\mathbf{u}$ and $\mathbf{x}$. Later, this splitting will be exploited to set up an offline-online decomposition in the computational steps, a feature of great importance for the efficiency of the reduced order method~\cite{Rozza:ARCME}.

\subsection{Full order approximation}
The full order approximation of equation~\eqref{eq:weakRBmu} with a Galerkin method reads: find $(\mathbf{u}^\mathcal{N},p^\mathcal{N})\in\mathbf{V}^\mathcal{N}\times Q^\mathcal{N}$
\begin{equation}
  \begin{cases}
    m(\mathbf{\mu};\mathbf{u}^\mathcal{N}(\mathbf{\mu}),\mathbf{v})+c(\mathbf{\mu};\mathbf{u}^\mathcal{N}(\mathbf{\mu}),\mathbf{u}^\mathcal{N}(\mathbf{\mu}),\mathbf{v})+a(\mathbf{\mu};\mathbf{u}^\mathcal{N}(\mathbf{\mu}),\mathbf{v}) &\\
    \qquad\qquad\qquad\qquad\qquad\qquad+b(\mathbf{\mu};\mathbf{v},p^\mathcal{N}(\mathbf{\mu}))=f(\mathbf{\mu};\mathbf{v})  &\pushright{\forall\mathbf{v}\in \mathbf{V}^\mathcal{N}} \\
    b(\mathbf{\mu};\mathbf{u}^\mathcal{N}(\mathbf{\mu}),q)=0 & \pushright{\forall q\in Q^\mathcal{N},}
  \end{cases}
  \label{eq:sem_ns}
\end{equation}
where $\mathbf{V}^\mathcal{N}\subset \mathbf{V}$ and $Q^\mathcal{N}\subset Q$ are the approximation spaces for the velocity and pressure, respectively.

For the construction of $\mathbf{V}^\mathcal{N}$ and $Q^\mathcal{N}$, we choose the Legendre Spectral Element Method (for details we refer for instance to~\cite{Fischer:02} or~\cite{CHQZ1,CHQZ2}) as implemented in the open source software \texttt{Nek5000}~\cite{nek5000}.
Velocity and pressure approximation spaces are constructed by Lagrangian interpolants of order 19 on Gauss-Lobatto-Legendre nodes distributed on a 12$\times$4 grid, for a total of 240$\times$80 collocation points in the two directions. As a result, the pressure approximation space has dimension $19200$, and the velocity approximation space has dimension $38400$.
For the time integration, the third order operator splitting method described in~\cite{Tombo} was used.

\section{A reduced order method applied to Navier--Stokes bifurcation problems}
\label{sec:ROM}
The aim of this section is to construct an approximation of equation~\eqref{eq:weakRBmu} with a Galerkin method  of dimension $N\ll\mathcal{N}$, so that the time to solution can be expected to decrease significantly. To this end, we recall the main features of a Reduced Basis Methods applied to equation~\eqref{eq:weakRBmu}.
The Reduced Basis formulation will be: find $(\mathbf{u}^N,p^N)\in\mathbf{V}^N\times Q^N$ such that
\begin{equation}
  \begin{cases}
    m(\mathbf{\mu};\mathbf{u}^N(\mathbf{\mu}),\mathbf{v})+c(\mathbf{\mu};\mathbf{u}^N(\mathbf{\mu}),\mathbf{u}^N(\mathbf{\mu}),\mathbf{v})+a(\mathbf{\mu};\mathbf{u}^N(\mathbf{\mu}),\mathbf{v}) &\\
    \qquad\qquad\qquad\qquad\qquad\qquad+b(\mathbf{\mu};\mathbf{v},p^N(\mathbf{\mu}))=f(\mathbf{\mu};\mathbf{v}) & \pushright{ \forall\mathbf{v}\in \mathbf{V}^N} \\
    b(\mathbf{\mu};\mathbf{u}^N(\mathbf{\mu}),q)=0 & \pushright{\forall q\in Q^N,}
  \end{cases}
  \label{eq:sem_ns}
\end{equation}
where $\mathbf{V}^N\subset \mathbf{V}^{\mathcal{N}}$ and $Q^N\subset Q^{\mathcal{N}}$ are the reduced basis spaces for the velocity and pressure respectively. In the next few paragraphs we discuss a strategy for the construction of $\mathbf{V}^N$ and $Q^N$.

%

\subsection{Sampling}
\label{sec:sampling}
The information needed for building a reduced order approximation is obtained by constructing a set of properly selected ``truth solutions'' $\{\mathbf{u}^{\mathcal{N}}(\mathbf{\mu}^i)\}_{i=1}^N\subset \mathbf{V}^{\mathcal{N}}$, where this term is customarily used to stress the fact that being the exact solutions $\mathbf{u}(\mathbf{\mu}^i)$ not accessible, their full-order approximations $\mathbf{u}^{\mathcal{N}}(\mathbf{\mu}^i)$ will be treated as if these were the exact solutions. The truth solutions are computed for a suitable sequence in the parameter space $\{\mathbf{\mu}^i\}_{i=1}^N\subset\mathcal{D}$. There exist different methods\footnote{We refer the interested reader to~\cite{Rozza:ARCME} for a brief overview on the history of the Reduced Basis Method and a review of many sampling techniques.} for identifying a sequence $\{\mathbf{\mu}^i\}\subset\mathcal{D}$, often based on a ``worst case'' criterion, which given an initial sequence $\{\mathbf{\mu}^i\}_{i=1}^k$, aims at searching for the parameter $\mu^{k+1}$ whose solution is the worst approximated one within the space spanned by the snapshots $\mathcal{S}_k=\{\mathbf{u}^{\mathcal{N}}(\mathbf{\mu}^i)\}_{i=1}^k$. Then, the new solution $\mathbf{u}^{\mathcal{N}}(\mathbf{\mu}^{k+1})$ is added to the snapshots space $\mathcal{S}_{k+1}=\mathcal{S}_k\cup\{\mathbf{u}^{\mathcal{N}}(\mathbf{\mu}^{k+1})\}$ and the algorithm is restarted until a suitable stopping criterion is satisfied.
 
A popular sampling method based on a worst case strategy is the \emph{Greedy Algorithm} (we refer to~\cite{GreedyNumerica} for a comprehensive review and to~\cite{Binev,Buffa} for some relevant convergence estimates).

In this work, we adopt the \emph{Centroidal Voronoi Tessellation (CVT)}, introduced in~\cite{Gunz} as sampling strategy for Reduced Order Modelling. The CVT builds a collection of subsets $\{\mathcal{V}_i\}_{i=1}^k$ of $\mathcal{D}$ (disjoint: $\mathcal{V}_i\cap\mathcal{V}_j=\emptyset$ if $i\ne j$, and where $\bigcup_{i=0}^k \mathcal{V}_i=\mathcal{D}$) such that, if we define the weight of each Voronoi region $\mathcal{V}_i$ by:
\begin{equation}
  \mathscr{W}_i=\int_{\mathcal{V}_i}\left\|\mathbf{u}(\mathbf{\nu})-\mathbf{u}(\xi_i)\right\|_{\mathbf{V}^{\mathcal{N}}}\left\|\mathbf{u}(\mathbf{\nu})-\Proj_{\mathbf{V}^N}\mathbf{u}(\mathbf{\nu})\right\|_{\mathbf{V}^{\mathcal{N}}}\de\mathbf{\nu},
        \label{eq:defCVTweight}
\end{equation}
where $\xi_i$ is the barycenter of the region $\mathcal{V}_i$:
        \begin{equation}
          \xi=\argmin_{\mathbf{\nu}\in\mathcal{U}}\int_{\mathcal{U}}\left\|\mathbf{u}(\mathbf{\mu})-\mathbf{u}(\mathbf{\nu})\right\|_{\mathbf{V}^{\mathcal{N}}}\left\|\mathbf{u}(\mathbf{\mu})-\Proj_{\mathbf{V}^N}\mathbf{u}(\mathbf{\mu})\right\|_{\mathbf{V}^{\mathcal{N}}}\de\mathbf{\mu}.
          \label{eq:Vmass}
        \end{equation}
        where $\Proj_{\mathbf{V}^N}:\mathbf{V}^{\mathcal{N}}\to \mathbf{V}^N$ is the orthogonal projector from $\mathbf{V}^{\mathcal{N}}$ to its subspace $\mathbf{V}^N$.
        The CVT of dimension $k$ of the parameter set is then defined by the generating points $\{\mathbf{\mu}^i\}_{i=1}^k$ that minimize the ``total weight'' functional:
\begin{equation}
        \mathscr{F}=\sum_{i=1}^k \mathscr{W}_i.
        \label{eq:totalW}
\end{equation}
In this sense the CVT is a best approximation sampling method, and is often combined with the POD defined in the next section to form the CVOD method~\cite{GunzCVOD}, quite popular in the ROM community.

For the sake of clarity, we report a scheme of the CVT algorithm adopted in this work in the Algorithm~\ref{algo:myCVT} box.

\begin{algorithm}
  \caption{The Centroidal Voronoi Tessellation sampling strategy used in Section~\ref{sec:stabilityreg}.}
  \label{algo:myCVT}
  \begin{algorithmic}[1]
   \Repeat
   \State{build the Delaunay triangulation for the points in the parameter space $\{\mathbf{\mu}^i\}\subset\mathcal{D}$;}
     \State{find the triangle $K_i\subset\mathcal{D}$ with the largest residual
     \begin{equation}
       \texttt{res}_{K_i}=\sum_{\mathbf{\mu}^j\in K_i}\|\mathbf{u}(\mathbf{\mu}^j)-\mathbf{u}^N(\mathbf{\mu}^j)\|_0;
     \end{equation}}
     \State{compute the barycenter $\mathbf{\mu^i}$ of $K_i$ and add it to the collection $\{\mathbf{\mu}^i\}$;}
     \State{compute the new snapshot $\mathbf{u}(\mathbf{\mu}^i)$;}
     \State{compute the corresponding basis function $\mathbf{\zeta}_i^{\mathrm{div}}$;}
     \State{update the tolerance:
       \begin{equation}
	 \texttt{tol}=\max_{i}\Big(\mathbf{u}(\mathbf{\mu}^i)-\sum_{j=1}^i(\mathbf{u}(\mathbf{\mu}^i),\mathbf{\zeta}_j^{\mathrm{div}})_0\mathbf{\zeta}_j^{\mathrm{div}}\Big);
         \label{eq:deftol}
       \end{equation}
     }
     \Until{\texttt{tol}$<$\texttt{threshold}.}
  \end{algorithmic}
\end{algorithm}

\subsection{Proper Orthogonal Decomposition}
\label{sec:POD}
When dealing with time-dependent problems, it is preferable not to add all the snapshots of a time-dependent run to the snapshots space $\mathcal{S}_k$, otherwise the Reduced Basis Space will likely be too large and many of the basis will be almost parallel, making the online projection phase an ill conditioned problem.
A common technique is to extract the ``most significant'' modes of a time sequence using a RB approximation in combination with a Proper Orthogonal Decomposition~\cite{Volkwein}.
There are many alternative ways to compute the POD modes of a sequence of snapshots $\{\mathbf{u}^{\mathcal{N}}_i\}\subseteq \mathbf{V}^{\mathcal{N}}$. Here we focus on the method based on the eigenvalues of the correlation matrix~\cite{Volkwein}. The entries of the correlation matrix $\mathbb{C}\in\mathbb{R}^{N\times N}$ are computed as
\begin{equation}
  \mathbb{C}_{ij}=(\mathbf{u}^{\mathcal{N}}_i,\mathbf{u}^{\mathcal{N}}_j)_{\mathbf{V}^{\mathcal{N}}},
  \label{eq:Cij}
\end{equation}
then, if $(\lambda_i,\psi_i)$ is a pair of eigenvalue and eigenvector of $\mathbb{C}$, each basis vector is computed as
\begin{equation}
  \mathbf{\zeta}_i=\sum_{k=1}^N \psi_{i,k} \mathbf{u}^{\mathcal{N}}_k
  \label{eq:PODz}
\end{equation}
where $\psi_{i,k}$ denotes the $k$-th component of the $i$-th eigenvalue.
The POD modes obtained are automatically orthogonal, but not normal in general. The eigenvalue $\lambda_i$ associated to each POD mode is related to the fraction of energy stored in the corresponding mode.

A remarkable property of the space generated by the POD bases is that it minimizes the projection error in the norm of the space $\mathbf{V}$.

The sampling algorithms presented in Section~\ref{sec:sampling} are frequently combined with a POD for parametrized time-dependent problems. If a Greedy sampling algorithm is used to select the parameters, and the POD is used to recover the most relevant time snapshots, the sampling procedure is called POD-Greedy and we refer to~\cite{HO08} and~\cite{CalcoloNRP} for details. 
To clarify the sampling process of time-dependent parametrized problems, we report in Algorithm~\ref{alg:PODGreedy} a possible implementation strategy of the POD-Greedy algorithm.

\begin{algorithm}
  \caption{A POD-CVT strategy for the sampling of parameter dependent evolution problems.}
  \label{alg:PODGreedy}
  \begin{algorithmic}[1]
   \Repeat
   \State{find $\mathbf{\mu}^{i+1}$ with the CVT Algorithm~\ref{algo:myCVT};}
   \State{compute a sequence of snapshots $\{\mathbf{u}_{i+1,k}\}_{k=0}^r$ for the time dependent problem related to the parameter $\mathbf{\mu}^i$;}
     \State{compute the $\ell$ POD modes of the sequence $\{\phi_{i+1,k}\}_{k=1}^\ell$ such that the retained energy is above a prescribed ratio;}
     \State{orthogonalize the time modes with respect to the previous basis sets: $\mathbf{\zeta}_{i+1,k}=\mathbf{\phi}_{i+1,k}-\Proj_{\mathcal{S}_i\cup_{j=1}^k\mathbf{\phi}_{i+1,j}}\mathbf{\phi}_{i+1,k}$;}
     \State{add the new basis functions $\{\mathbf{\zeta}_{i+1,k}\}_{k=1}^\ell$ to the basis set: $\mathcal{S}_{i+1}=\mathcal{S}_i\cup\{\mathbf{\zeta}_{i+1,k}\}$;}
   \Until{$\max_{\mu\in\mathit{\Sigma}}\Delta(\mu) < \texttt{tol}$.}
  \end{algorithmic}
\end{algorithm}

Lastly, we remark that when sampling time dependent problems in view of a POD application, it is important to make sure that the sampling rate is sufficiently high so that the desired time harmonics are well resolved. 

\subsection{Reduced order formulation of the parametrized Navier-Stokes equations}
\label{sec:reducedNS}
We remark that as for the truth spaces, the reduced basis spaces should be chosen so that three fundamental properties are fulfilled. First, it is important that good stability properties are verified, that is, the approximation spaces must lead to a well-posed problem (\emph{approximation stability}). Secondly, the reduced order spaces should guarantee good approximation properties for all the parameters on a given interval. Third, the reduced order spaces should have the lowest possible dimension while mantaining the required \emph{approximation properties}.
In the following we will discuss how these properties can be satisfied by a careful choice of the spaces $\mathbf{V}^N$ and $Q^N$, respectively the velocity and pressure reduced basis spaces.

A basis for the spaces $\mathbf{V}^N$ and $Q^N$ can be computed through the techniques discussed in Section~\ref{sec:POD}, but some additional care is required to ensure the approximation stability of the reduced basis spaces. In general, the basis obtained as described in Section~\ref{sec:POD} will not fulfill the LBB inf-sup condition:
\begin{equation}
        \beta^N\equiv\inf_{q\in Q^N}\sup_{\mathbf{v}\in\mathbf{V}^N}\frac{b(\mathbf{\mu};q,\mathbf{v})}{\|q\|_{Q^N}\|\mathbf{v}\|_{\mathbf{V}^N}}>0.
  \label{eq:Brezziinfsup}
\end{equation}

Two popular possibilities to recover the inf-sup control are the \emph{supremizer enrichment} (see for instance~\cite{RozzaVeroy} for an introduction, and~\cite{RHM} for an analysis of the supremizer stabilization), or the \emph{Petrov-Galerkin stabilization} (we refer to~\cite{Rovas} for general nonsymmetric problems, and to~\cite{Dahmen} and~\cite{AbdullePG} for Stokes problems).

In this work, we adopt the \emph{Piola transformation} technique, introduced in the ROM community in~\cite{Lovgren} which consists in a suitable preprocessing of the velocity snapshots $\{\mathbf{u}^\mathcal{N}_i\}$ that returns a set of divergence-free velocity basis functions $\{\mathbf{\zeta}_i\}$ for each value of the geometric parameter $\mu$. The divergence-free basis set allows to cancel out the pressure term from the momentum equation, thus removing any stability issue for mixed problems. Pressure can then be recovered solving a Poisson problem online:
\begin{equation}
        \Delta p^N(\mathbf{\mu}^i)=-\ddiv\left(\mathbf{u}^N(\mathbf{\mu}^i)\cdot\nabla\mathbf{u}^N(\mathbf{\mu}^i)\right).
        \label{eq:prossurePoisson}
\end{equation}
We refer for example to~\cite{Caiazzo} for an analysis of velocity-pressure reduced order models.

The Piola transformation is based on the following consideration. Each snapshot $\mathbf{u}^{\mathcal{N}}(\mathbf{\mu}^i)$ is weakly divergence-free in the original domain $\mathit{\Omega}(\mu^i)$, but this is not true when the snapshot is pulled back to the reference domain $\widehat{\mathit{\Omega}}$. Indeed, the parametrized formulation imposes:
\begin{equation}
  \int_{\mathit{\Omega}(\mu^i)}q\ddiv\mathbf{u}(\mathbf{\mu}^i) \de\mathbf{x}=0\qquad\forall q\in Q^N,
  \label{eq:HOdiff}
\end{equation}
that in coordinates reads:
\begin{equation}
  \begin{split}
    \int_{\mathit{\Omega}(\mu^i)}q\Big(\sum_{j=1}^d\frac{\partial\mathbf{u}_j(\mathbf{\mu}^i)}{\partial x_j}\Big)\de\mathbf{x}=
    \int_{\widehat{\mathit{\Omega}}}q\Big(\sum_{j=1}^d\sum_{k=1}^dG_{jk}(\mu^i)\frac{\partial\mathbf{u}_j(\mathbf{\mu}^i)}{\partial \widehat{x}_k}\Big)J^{\mathrm{aff}}(\mu^i)\de\widehat{\mathbf{x}}.
  \end{split}
  \label{eq:HOdiffcoordinates}
\end{equation}
It can be concluded that the snapshots do not cancel out the standard divergence on the reference domain:
\begin{equation}
        \int_{\widehat{\mathit{\Omega}}}q\sum_{j=1}^d\frac{\partial \mathbf{u}_j(\mathbf{\mu}^i)}{\partial \widehat{x}_j}\de\widehat{\mathbf{x}}\qquad\forall q\in Q^N
  \label{eq:originaldiv}
\end{equation}
but instead the pushed forward, or ``stretched'' divergence~\eqref{eq:HOdiffcoordinates}.


The algorithm of the Piola transformation can be summarized as follows.
\begin{enumerate}
  \item the snapshots $\{\mathbf{u}^\mathcal{N}(\mathcal{\mu}^i)\}$ are mapped to the reference domain $\widehat{\mathit{\Omega}}$ by terms of the inverse Piola transformation $\mathcal{P}^{-1}(\mu)$, so that after the pull-back the snapshots $\{\mathcal{P}^{-1}(\mu^i)\mathbf{u}^\mathcal{N}(\mu^i)\}$ are divergence-free on the reference domain.
  \item the RB functions $\{\widehat{\mathbf{\zeta}}\}$ are computed on the reference domain, for instance via a POD procedure starting from the pulled-back snapshots.
  \item the mass, stiffness and convection matrices are assembled on the reference domain, for instance for the mass and stiffness matrices:
\begin{equation}
        \begin{aligned}
                &\widehat{M}_{ij}=\int_{\widehat{\mathit{\Omega}}}\widehat{\mathbf{\zeta}}_i\cdot\widehat{\mathbf{\zeta}}_j\de\widehat{\mathbf{x}}
                \qquad\widehat{K}_{ij}=\int_{\widehat{\mathit{\Omega}}}\widehat{\nabla}\widehat{\mathbf{\zeta}}_i\cdot\widehat{\nabla}\widehat{\mathbf{\zeta}}_j\de\widehat{\mathbf{x}}\\
        \end{aligned}
        \label{eq:defmatrices}
\end{equation}
note that it is not required to assemble any pressure-related matrix.

\item During the online phase, the variational forms~\eqref{eq:bilinearforms},~\eqref{eq:formC},~\eqref{eq:defF} should be mapped on the parametrized divergence-free basis $\{\mathbf{\zeta}_i\}$. Such basis are obtained through a Piola transformation of the original basis:
\begin{equation}
        \begin{sistema}
                \mathbf{\zeta}_{i,1}=\mathcal{P}_{1}(\mu)\widehat{\mathbf{\zeta}}=G_{11}(\mu)\widehat{\mathbf{\zeta}}_{i,1}+G_{12}(\mu)\widehat{\mathbf{\zeta}}_{i,2}\\
                \mathbf{\zeta}_{i,2}=\mathcal{P}_{2}(\mu)\widehat{\mathbf{\zeta}}=G_{21}(\mu)\widehat{\mathbf{\zeta}}_{i,1}+G_{22}(\mu)\widehat{\mathbf{\zeta}}_{i,2}.\\
        \end{sistema}
        \label{eq:defparamLerayP}
\end{equation}
During the parametrized simulation, it is not necessary to explicitly build the new basis $\{\mathbf{\zeta}_i\}$, but only to evaluate its effect on the matrices~\eqref{eq:defmatrices} (i. e. on the coeffcients). The expression of the new matrices can be obtained from the reference matrices~\eqref{eq:defmatrices} by a change of variables under the integral sign, as done in equation~\eqref{eq:HOdiffcoordinates}.

\end{enumerate}

In this work, we explore the possibilities offered by the Piola transformation to build a divergence-free basis on each parametrized domain.

\section{Numerical detection of bifurcation points}
\label{sec:bifurcation_detection}

In this section we define two types of bifurcations: the steady bifurcation and the Hopf bifurcation. Then we recall briefly some numerical techniques for the detection of these two types of bifurcations, and we discuss their applicability in a reduced order framework.

\subsection{Steady bifurcation}
\label{sec:steady}
Intuitively, for a fixed value of the domain length $A$, we say that $\Gr^*$ is a \emph{bifurcation point} for~\eqref{eq:weakRBmu} if the equation~\eqref{eq:weakRBmu} admits two distinct solutions $(\mathbf{u}^1,p^1)$ and $(\mathbf{u}^2,p^2)$ for $\Gr\ge\Gr^*$, but only one solution for $\Gr<\Gr^*$.
This qualitative observation can be made formal~\cite{Prodi} by defining $\Gr^*$ a \emph{bifurcation point} if there exist two distinct sequences of parameter-solution couples $(\Gr_m^1,(\mathbf{u}_m^1,p_m^1))$ and $(\Gr_n^2,(\mathbf{u}_n^2,p_n^2))$ such that:
\begin{itemize}
  \item[i.] $(\Gr_m^1,(\mathbf{u}_m^1,p_m^1))$ and $(\Gr_n^2,(\mathbf{u}_n^2,p_n^2))$ are solutions of the equation~\eqref{eq:weakRBmu};
  \item[ii.] if $(\mathbf{u}^*,p^*)$ is the unique solution of equation~\eqref{eq:weakRBmu} for the parameters $(\Gr^*,A)$, then:
    \begin{equation}
      (\Gr_m^1,(\mathbf{u}_m^1,p_m^1))\to (\Gr^*,(\mathbf{u}^*,p^*)) \quad \text{and} \quad (\Gr_n^2,(\mathbf{u}_n^2,p_n^2))\to(\Gr^*,(\mathbf{u}^*,p^*))
      \label{eq:def_bif}
    \end{equation}
\end{itemize}
as $m,n\to\infty$.

To identify numerically a steady state bifurcation point, we follow a technique similiar to Lemma~4 in~\cite{Galdi:NSanalysis} that we briefly summarize.
We introduce the linearized advection operator $\mathbf{\mathcal{T}}(\mathbf{u}^*):\mathbf{V}\mapsto\mathbf{V}$, obtained by taking the Fr\'echet derivative of the convection term $\mathbf{u}\cdot\nabla\mathbf{u}$ about a base solution $\mathbf{u}^*$:
\begin{equation}
  \mathbf{\mathcal{T}}(\mathbf{u}^*)[\mathbf{v}]\equiv\mathbf{u}^*\cdot\nabla\mathbf{v}+\mathbf{v}\cdot\nabla\mathbf{u}^*.
  \label{eq:Tdef}
\end{equation}
According to~\cite{Galdi:NSanalysis}, if $(\mu^*,\mathbf{u}^*)$ is a bifurcation point ($\mu$ is a parameter), the equation
\begin{equation}
  \mu^*\mathbf{v}+\mathbf{\mathcal{T}}(\mathbf{v}^*(\mu^*))[\mathbf{v}]=0
  \label{eq:galdieq}
\end{equation}
has at least one nonzero solution $\mathbf{v}$.

From a reduced order modelling perspective, following~\cite{cliffe} we search for a change of sign of the eigenvalues of the matrix $T(\mathbf{u}^*)$, defined as
\begin{equation}
  \begin{aligned}
    T_{ij}(\mathbf{u}^*)&=\mathbf{\mathcal{T}}((\mathbf{u}^*,\mathbf{\zeta}_j)_{L^2}\mathbf{\zeta}_j)[\mathbf{\zeta}_i]\\
    &=\sum_{k=1}^N(\mathbf{\zeta}_i,\mathbf{\zeta}_k\cdot\nabla\mathbf{\zeta}_j)_{L^2}\mathbf{U}^k_N+\sum_{k=1}^N(\mathbf{\zeta}_i,\mathbf{\zeta}_j\cdot\nabla\mathbf{\zeta}_k)_{L^2}\mathbf{U}^k_N,
\end{aligned}
  \label{eq:defmatrixT}
\end{equation}
where $(\cdot,\cdot)_{L^2}$ denotes the inner product in $L^2(\Om)$.

One important remark has to be made regarding the choice of the snapshots from which the reduced basis shall be selected.
To detect the null eigenvalue of the matrix $T(\mathbf{u}^*)$ of equation~\eqref{eq:defmatrixT}, it is important that $\mathbf{\zeta}_i\in\ker\mathbf{\mathcal{T}}( ( \mathbf{u}^*,\mathbf{\zeta}_j )_{L^2},\mathbf{\zeta}_j)$ for at least one couple $(i,j)$. Otherwise, all the eigenvalues of $T(\mathbf{u}^*)$ will be nonzero for all the values of $\mathbf{\mu}$, and no bifurcation will be detected.

In general it may be difficult to compute a basis function for the kernel of $\mathbf{\mathcal{T}}$ and add it to the space $\mathbf{V}^N$, but from numerical experiments we can suggest an alternative practical way to fulfil this condition.

We claim that it is sufficient to select the reduced basis functions from snapshots coming both from the bifurcating branch and from the branch of solutions before the bifurcation, and the POD procedure among the other will select from the snapshots an approximation for the kernel of $\mathbf{\mathcal{T}}$.

In the following, we present some numerical evidence to support such claim.
A first experiment consists in the reduced order computation of a steady bifurcation point for the case $A=4$ where the reduced basis function come from snapshots coming both from the single-vortex branch and the two-vortex branch.

Then, we try to detect again the bifurcation point with the reduced order model but this time we perform the POD on the snapshots coming from the two-vortex branch only.

\begin{figure}[htbp]
  \centering
  \includegraphics[width=0.49\textwidth]{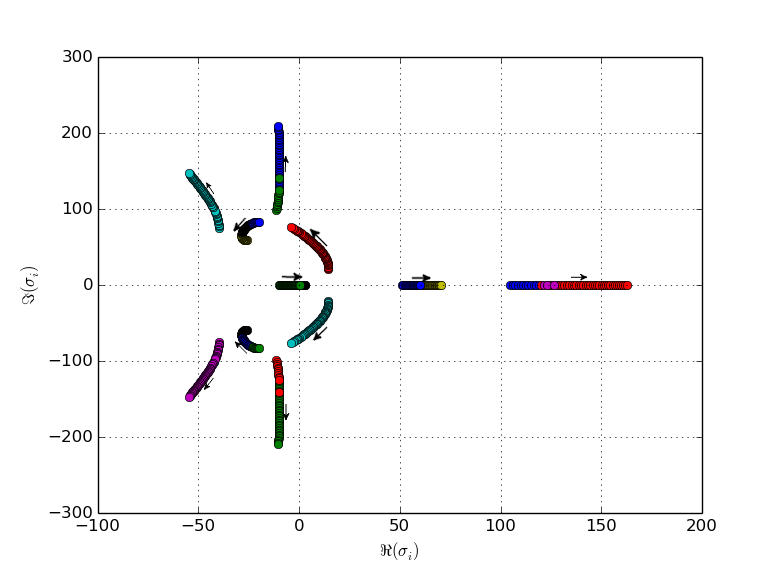}
  \includegraphics[width=0.49\textwidth]{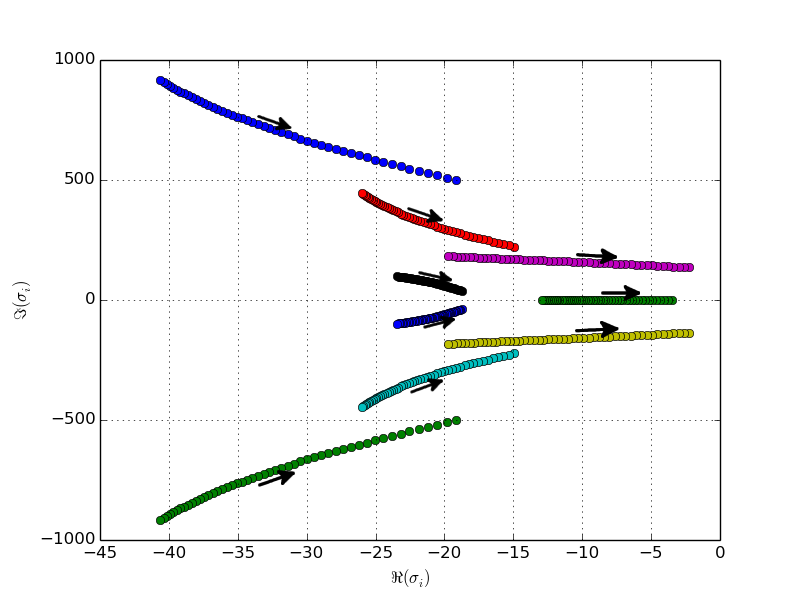}
  \caption{Reduced Model Eigenvalues of the tangent advection operator $\mathbf{\mathcal{T}}$ at the bifurcation point for $A=4$, for $\Gr\in[80\cdot10^3,110\cdot10^3]$. 
  On the left the tangent operator is evaluated on a basis set with one, two and three rolls, on the right a set with only two and three roll flows. The colors have the only purpose of separating visually the different eigenvalues.}
  \label{fig:eigsLin}
\end{figure}

The path of the eigenvalues of the operator $\mathbf{\mathcal{T}}$ in a neighbourhood of the bifurcation point are shown in Figure~\ref{fig:eigsLin}. From the picture it can be seen that while in the first case a real negative eigenvalue passes through the origin of the plane to become real positive, in the second test case all the eigenvalues keep a negative real part, thus no bifurcation is detected.

\begin{figure}[htbp]
  \centering
  \includegraphics[width=0.55\textwidth]{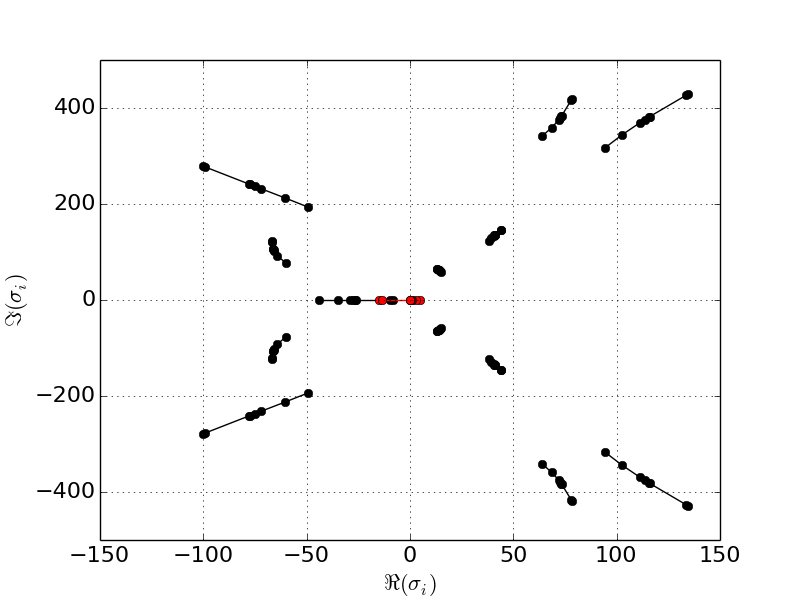}
  \caption{Evolution in the complex plane of the Reduced Model eigenvalues of $\mathbf{\mathcal{T}}$ with the Grashof number  in a neighborhood of a steady bifurcation point. Highlighted in red is the path of the critical eigenvalue, responsible for the bifurcation.}
  \label{fig:alleigs}
\end{figure}

As a further example, we report in Figure~\ref{fig:alleigs} the path of the eigenvalues in a neighborhood of a steady bifurcation point. Highlighted in red is the path of the eigenvalue changing its sign as the Grashof number is increased. 

We remark that since $\mathrm{T}_{ij}$ has a small dimension, all the eigenvalues can be computed at a reasonable expense, for instance $O(N^3)$ using QR iterations~\cite{Golub}. 




We do not make any claim regarding the relation of the RB eigenvalues and the full order eigenvalues. In fact, the path of the RB eigenvalues in the complex plane as $\mu$ is varied may be significantly different from the path followed by the full order eigenvalues. This however is not an issue, since we are not interested in reproducing correctly the path of the full order spectra\footnote{Which is usually dependent on the discretization parameter such as grid size.}, but only in correctly computing the value of $\mu$ for which an eigenvalue crosses the imaginary axis.

Recently, a POD-based eigenvalue approach has been introduced for the detection of mechanical vibration in the automotive industry~\cite{Mehrmann}, showing reliable results.

In any case, we remark that the bifurcating behaviour of the RB solution can be seen only if the samples are coming from two branches, i.e. the bifurcation has appeared also in the full order setting.

\subsection{Hopf bifurcation}

The Hopf bifurcation does not regard the branching of steady state solutions, instead it identifies a point in the parameter space for which a steady state solution becomes time-dependent.
This implies that at a Hopf bifurcation point a random perturbation from a steady state does not damp off in time, but after a transient phase it leads to a new, time-dependent solution. 

For the detection of Hopf bifurcation points, we write the time-dependent Navier-Stokes equations seeking for a solution in the form of a superposition of a steady state solution  $\mathbf{u}^*(\mathbf{x})$ and a small perturbation $\mathbf{u'}(\mathbf{x})\mathrm{e}^{\sigma t}$, $\sigma\in\mathbb{C}$. Neglecting the second order terms, we obtain a linearized equation that rules the time evolution of the perturbations:
\begin{equation}
  \mathbf{u}^*\cdot\nabla\mathbf{u'}+\mathbf{u'}\cdot\nabla\mathbf{u}^*-\Delta\mathbf{u'}=-\sigma \mathbf{u'}
  \label{eq:NSperturb}
\end{equation}
that can be seen as an eigenvalue problem for the linearized Navier-Stokes operator $\mathbf{\mathcal{L}}(\mathbf{u}^*):\mathbf{V}\mapsto\mathbf{V}$:
\begin{equation}
        \mathbf{\mathcal{L}}(\mathbf{u}^*)[\mathbf{u'}]=-\sigma\mathbf{u'}.
        \label{eq:defLNS}
\end{equation}
If equation~\eqref{eq:defLNS} admits an eigenvalue $\sigma^*$ such that $\Re\sigma^*>0$, the corresponding perturbation will grow in time, and if $\Im\sigma^*>0$, an oscillatory solution has to be expected, at least until the nonlinear term $\mathbf{u'}\cdot\nabla\mathbf{u'}$ remains sufficiently small.

Applying this technique to the reduced order model, the Hopf bifurcation is detected when the matrix $L$ associated to the linearized Navier-Stokes operator $\mathbf{\mathcal{L}}(\mathbf{u}^*)$ admits an eigenvalue satisfying the conditions above.
Explicitly, the matrix $L$ has the form
\begin{equation}
  L_{ij}=\sum_{k=1}^{N_{\mathbf{u}}}(\mathbf{\zeta}_i,\mathbf{\zeta}_k\cdot\nabla\mathbf{\zeta}_j)_{L^2}\mathbf{U}^k_N+\sum_{k=1}^{N_{\mathbf{u}}}(\mathbf{\zeta}_i,\mathbf{\zeta}_j\cdot\nabla\mathbf{\zeta}_k)_{L^2}\mathbf{U}^k_N+(\nabla\mathbf{\zeta}_i,\nabla\mathbf{\zeta}_j)_{L^2}.
  \label{eq:matrixL}
\end{equation}

As for the operator $\mathcal{T}$ discussed in Section~\ref{sec:steady}, we remark that:
\begin{itemize}
  \item the matrix $L$ of equation~\eqref{eq:matrixL} associated to the operator $\mathbf{\mathcal{L}}(\mathbf{u}^*)$ is sufficiently small that no Krylov methods are needed to compute its eigenvalues;
  \item it is important to select snapshots both from the steady state and the time-dependent solution branches, otherwise the Hopf bifurcation may not be detected in the reduced order setting.
\end{itemize}

\section{Numerical experiments}
\label{sec:results}

We now test the efficiency of the reduced basis method for flow bifurcation problems on the test case described in Section~\ref{sec:appfluid}. Many experimental and numerical references exist on this benchmark, among which we will refer to the results in~\cite{Gelfgat:Ref11}.

In this benchmark the parameters are the domain length $A\in[2,10]$ and the Grashof number $\Gr\in[50\cdot10^3,1\cdot10^6]$. 
In the considered parameter domain there exist regions admitting unique steady solutions, multiple steady solutions, and unsteady solutions, consequently this represents a good problem to test the proposed numerical method.

\begin{figure}[tbp]
\centering
\includegraphics[width=0.2\textwidth]{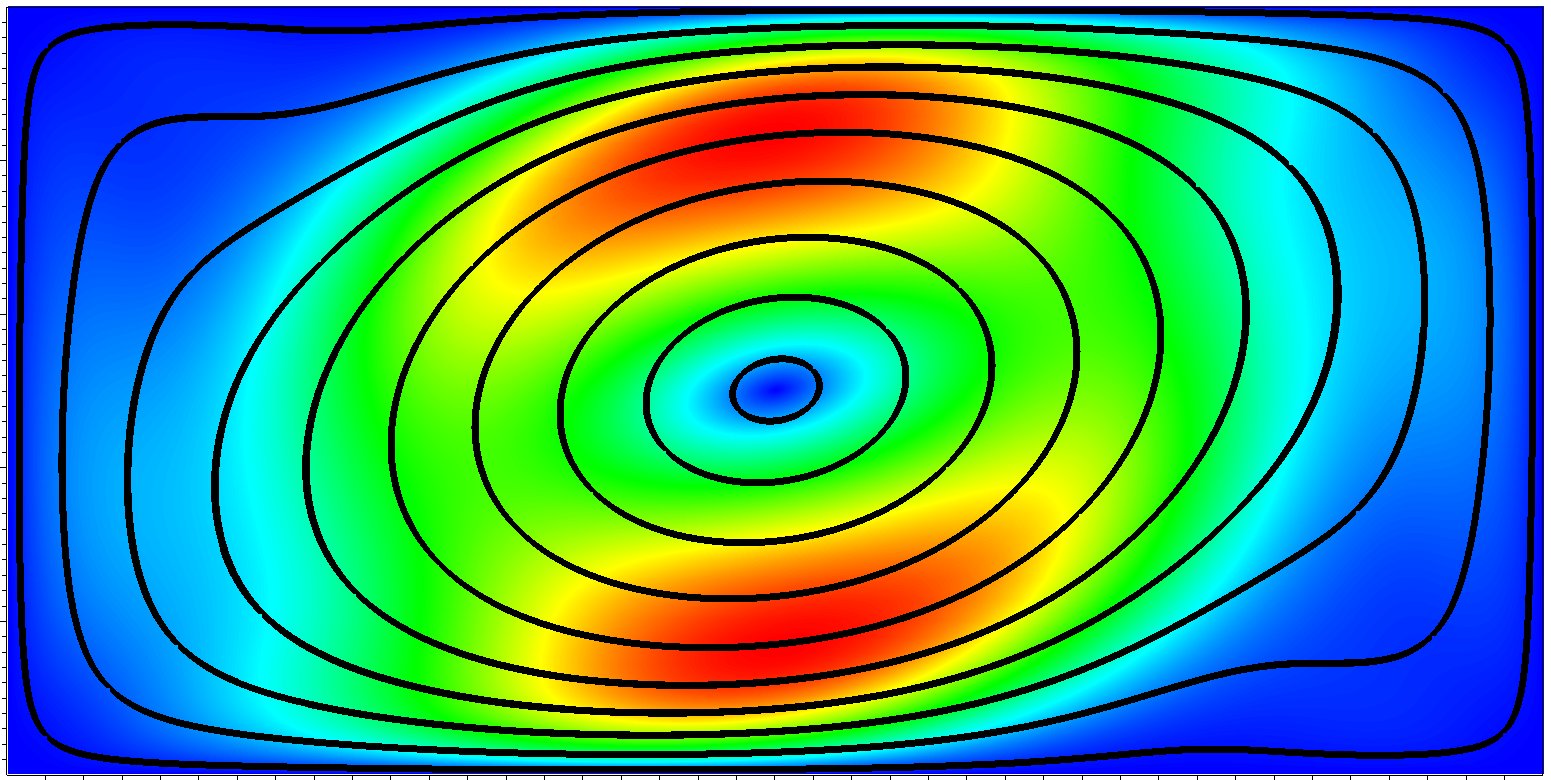}
\includegraphics[width=0.337\textwidth]{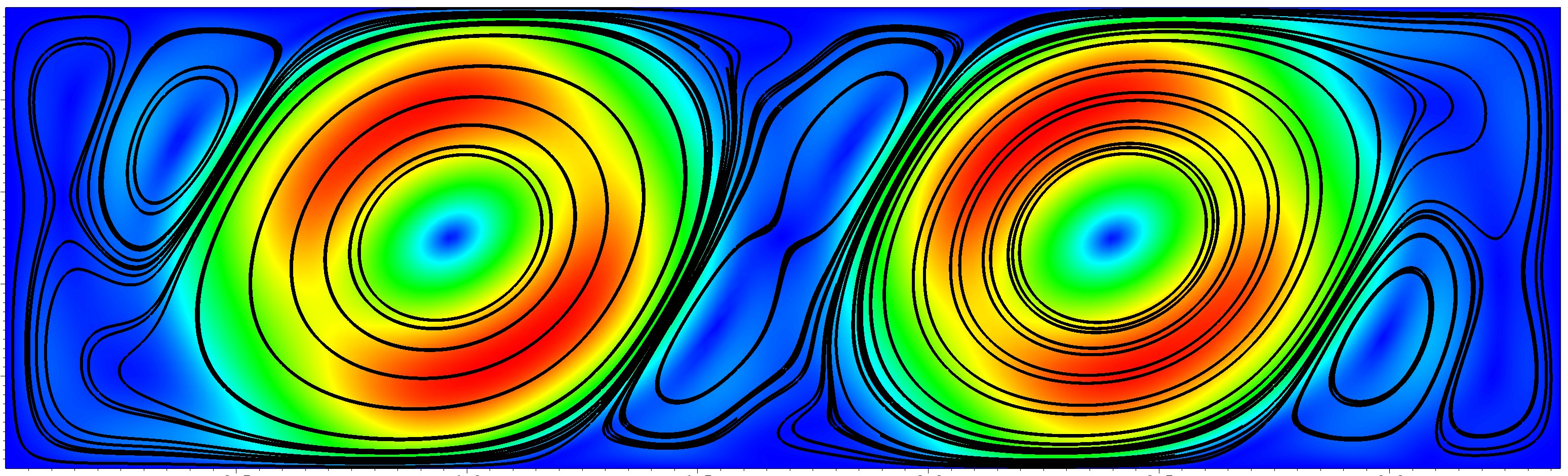}
\includegraphics[width=0.552\textwidth]{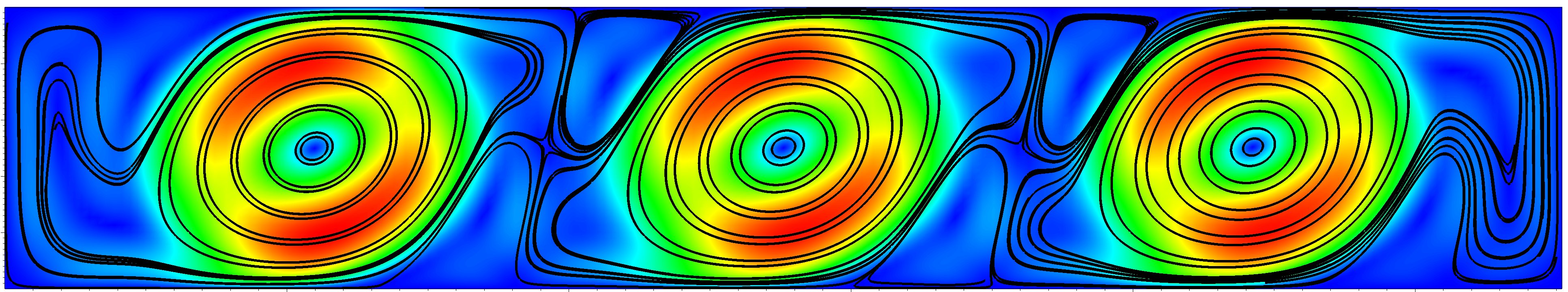}
\includegraphics[width=0.836\textwidth]{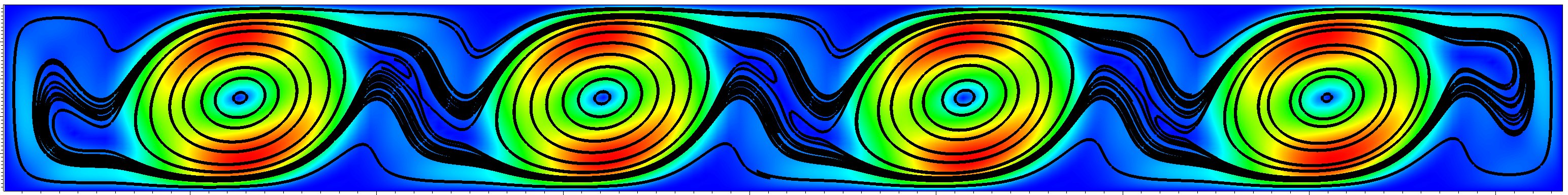}
\includegraphics[width=\textwidth]{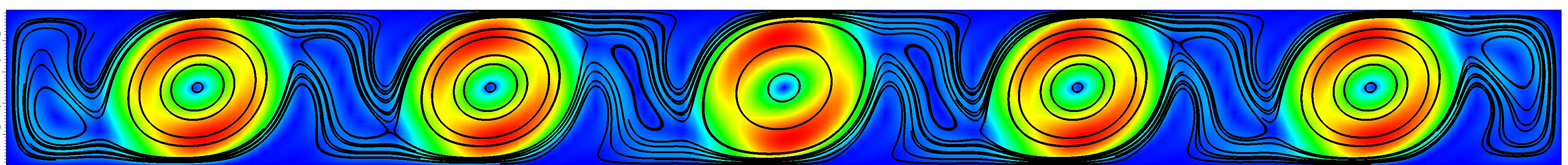}
\caption{Some snapshots for the GAMM benchmark. From top to bottom, and from left to right: $(A=2,\Gr=50\cdot10^3)$; $(A=3.37,\Gr=264.9\cdot10^3)$; $(A=5.52,\Gr=132.1\cdot10^3)$; $(A=8.36,\Gr=50.77\cdot10^3)$; $(A=10,\Gr=100\cdot10^3)$. All snapshots are steady state, except the one with three rolls.}
\label{fig:snapsGAMM}
\end{figure}

Some representative solutions for this benchmark are shown in Figure~\ref{fig:snapsGAMM}. As it can be seen in the visualizations, the solutions are very heterogeneous as the Grashof number and aspect ratio are varied within $\mathcal{D}$.

We choose as reference domain the rectangle with $A=2$, $\widehat{\mathit{\Omega}}=\mathit{\Omega}(2)$ so the affine transformation is given by:
\begin{equation}
  \begin{pmatrix}
    x \\
   y
  \end{pmatrix} =
  \mathcal{T}^{\mathrm{aff}}(\widehat{\mathbf{x}},A)=
  \begin{pmatrix}
  0 \\
  0
  \end{pmatrix}
  +
  \begin{bmatrix}
   \frac{A}{2} & 0 \\
  0 & 1  \\
  \end{bmatrix}
  \begin{pmatrix}
          \widehat{x} \\
  \widehat{y}
  \end{pmatrix}.
  \label{eq:mytransform}
\end{equation}
For this simple geometry, the Piola transformation $\mathcal{P}(A,\mathbf{x})$, is simply the inverse of $G(A)$.


For the time evolution of the online system we chose a third order semi-implicit Backward Difference Formula (BDF3).

\subsection{Single-parameter test}
\label{sec:single_parameter_test}
In the first test case the geometric aspect ratio is fixed to $A=4$, and the Grashof number is the only active parameter. 
According to~\cite{Gelfgat:Ref11}, for this aspect ratio there are three steady solutions up to $\Gr=120\cdot10^{3}$, after which a Hopf bifurcation occurs and the flow becomes unsteady. The three branches of steady solutions are characterized by a single roll flow up to $\Gr=25\cdot10^3$, two rolls for $\Gr\le100\cdot10^3$, and three rolls from this point up to the Hopf bifurcation at $\Gr=120\cdot10^3$.

The tolerance for the sampling procedure is $10^{-4}$, and the sampling procedure provided 13 snapshots for $\Gr\in[40\cdot10^3,1\cdot10^6]$, of which 7 are steady state solutions, and the remaining 6 are computed by a POD of 2 time periodic solutions, with an $L^2$ energy threshold fixed to 99.9\%.
To test the proposed methodology (CVT-POD, SEM, Piola transformation) in the usual RB framework we try to rebuild a known bifurcation diagram with the reduced model, with the aim to reconstruct the different solution branches.



\begin{figure}[htbp]
  \centering
  \includegraphics[width=\textwidth]{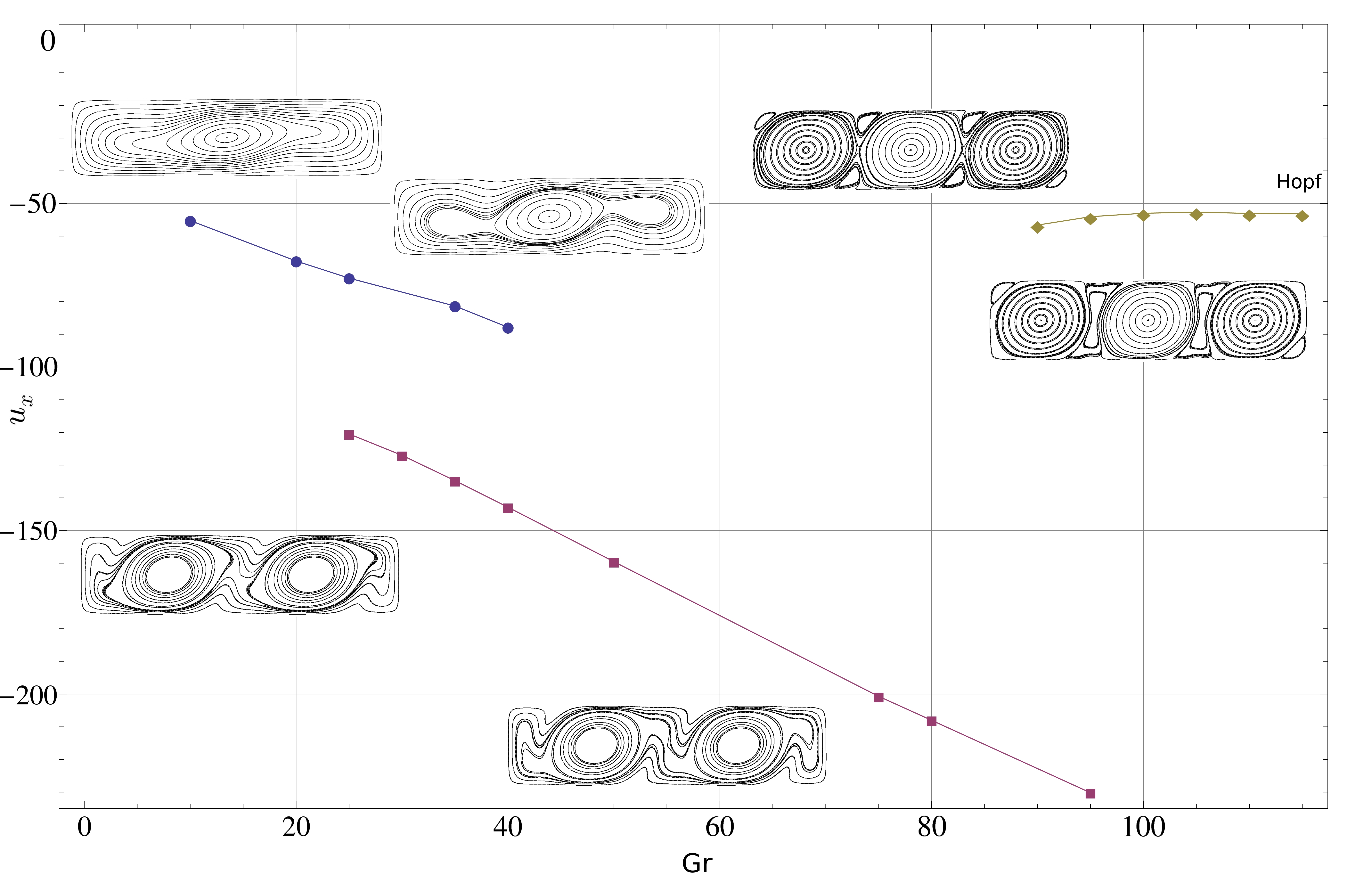}
  \caption{Bifurcation diagram for an aspect ratio of 4. The graph shows the horizontal velocity at the point $(0.7,0.7)$ as a function of the Grashof number. The three lines are associated to the solutions with 1, 2, and 3 rolls, and the streamlines of some representative solutions are plotted for clarity.}
  \label{fig:A4diagram}
\end{figure}

To build the online bifurcation diagram, we run the continuation method until the Hopf bifurcation is reached. To check for the presence of hysteresis, the continuation method is run backwards, i.e. decreasing the Grashof number from $\Gr=120\cdot10^3$.
The bifurcation diagram is shown in Figure~\ref{fig:A4diagram}, where the evolution of the horizontal velocity at a fixed point inside the domain is plot as a function of the Grashof number. The three branches can clearly be identified, and some hysteresis is present, especially at the transition between one and two rolls flows. For clarity reasons, we also plot the streamlines of some representative solutions.

\begin{figure}[tbp]
\centering
\includegraphics[width=0.4\textwidth]{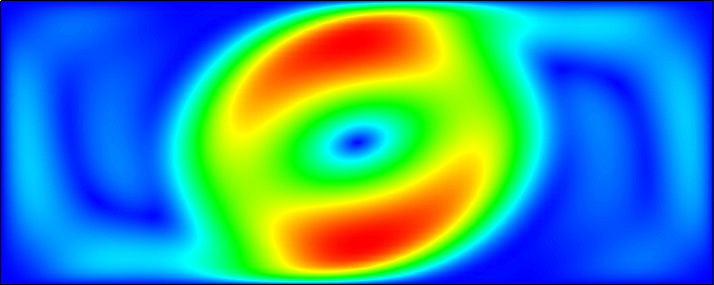}
\includegraphics[width=0.4\textwidth]{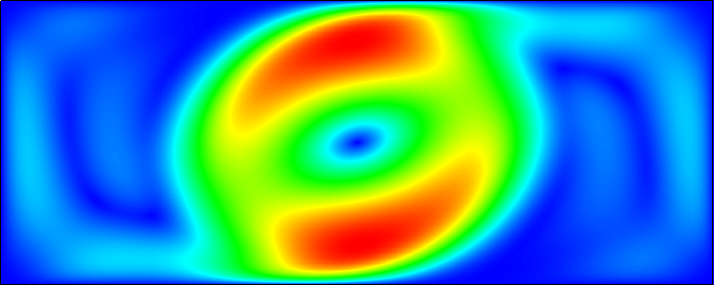} \\
\includegraphics[width=0.5\textwidth]{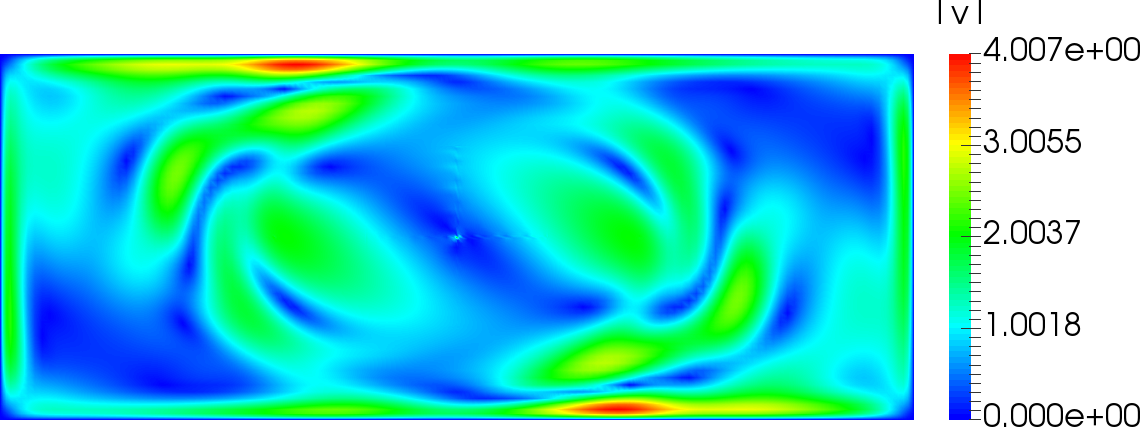} \hfill
\caption{A comparison of the full order (top row, left) and the reduced order (top row, right) snapshots in the time-periodic case with parameters $\Gr=270\cdot10^3$, $A=2.5$.}
\label{fig:visual_cmp}
\end{figure}

A further visual check of the reduced basis methods is available in Figure~\ref{fig:visual_cmp}, where a comparison of the full order and the reduced order snapshots in the time-periodic case with parameters $\Gr=270\cdot10^3$, $A=2.5$, as well as the relative $L^2$ distance between the two solutions is shown.

The online procedure required on average slightly less than 600 seconds, while the full order run would have taken about 24 cpu-hours for each solution on a IBM iDataPlex DX360M3 cluster, for a total of 576 cpu-hours. 

\subsection{Two-parameter tests}
\label{sec:stabilityreg}
After the fixed-geometry test, we consider the more general case with parametrized geometry. The parameter space $\mathcal{D}$ is now two dimensional, in particular we  choose $\Gr\in[50\cdot10^3,1\cdot10^6]$ and $A\in[2,6]$. 


The parametrized geometry requires to enable the Piola transformation during the online phase; this was not the case in Section~\ref{sec:single_parameter_test} since there the geometry was fixed.

\begin{figure}[htbp]
  \centering
  \includegraphics[width=0.75\textwidth]{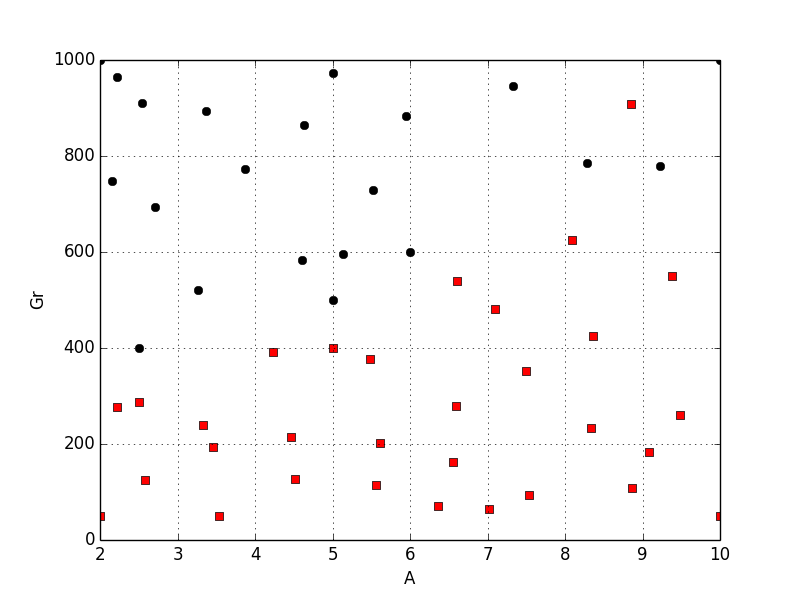}
  \caption{Parameters selected by the sampling Algorithm~\ref{alg:PODGreedy} of Section~\ref{sec:sampling}. The red square mark denotes the steady state snapshots, the black circle marks the time-dependent snapshots. The Grashof number is expressed in thousands.}
  \label{fig:mysnaps}
\end{figure}

In this test, the sampling algorithm provided a total of 108 basis functions, arising from 30 steady and 21 unsteady snapshots. Almost all the unsteady snapshots required 3 POD modes to store the 99.9\% of the energy, except for a few cases with 2 modes.
We report in Figure~\ref{fig:mysnaps} the position of the snapshots on the parameter plane $(A,\Gr)$, as computed by the POD-CVT algorithm.

\begin{figure}[htbp]
        \centering
        \includegraphics[width=0.75\textwidth]{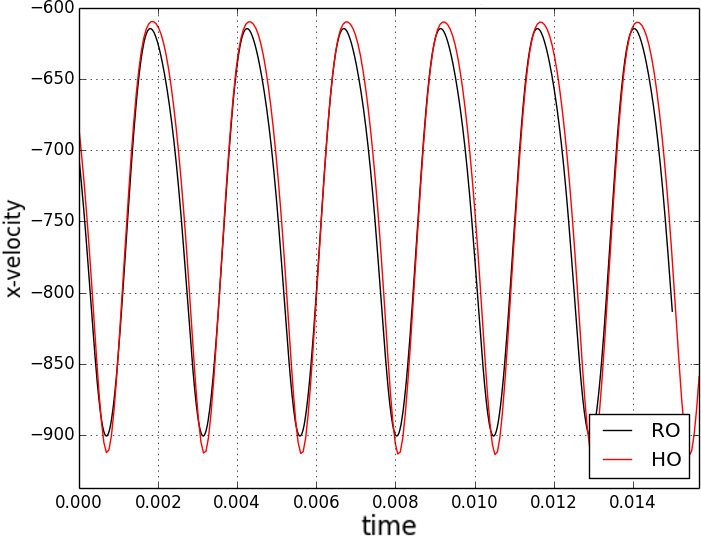}
        \caption{Comparison of the horizontal velocity at a the point $(0.7,0.7)$ vs time for the high order (in red, HO) and reduced order (in black, RO) simulations. The parameters are $\Gr=963791$, $A=2.22$, and the resulting flow has a single roll.}
        \label{fig:dinamicsR1}
\end{figure}

The first validation for the two parameter case consists in the reproduction of a time-periodic regime. We choose the point in the parameter space with coordinates $A=2.22$, $\Gr=963791$, and we run an online simulation with the reduced basis space $\mathbf{V}^N$ generated by the full 108 solutions. The online results are compared with the full-order data in Figure~\ref{fig:dinamicsR1}. The main features of the periodic regime are correctly reproduced, with an error in the estimated frequency of $0.51\%$. The computational savings for the single run are relevant:
\[
  \frac{\text{total offline and online computation time}}{\text{equivalent fully offline computation time}}\simeq\frac{5 \text{ cpu minutes}}{24 \text{ cpu hours}}=0.35\%.
        \label{eq:time_dependent_savings}
\]

The next test aims at the determination of steady and Hopf bifurcation points in the parameter plane, together with the frequency at the onset of the oscillations (i.e. immediately after the Hopf bifurcation).
The technique used to detect all the bifurcation points is the following:
\begin{description}
  \item[step-1:] fix an initial value of the geometrical parameter $A$;
  \item[step-2:] run the continuation algorithm for the fixed value of $A$;
  \item[step-3:] increase $A$ and go to step-2; if $A$ is in the boundary of the parameter space, change trial and test space and go to step-1.
\end{description}
To clarify the meaning of step-3, we explain how the procedure has been applied to detect the bifurcation points of the one-vortex flow. In a first run, the three steps have been carried out using a trial space obtained from the POD selection of the one-vortex snapshots only, and a test space obtained only from the two-vortex snapshots. After the bifurcation points with such approximation spaces have been computed for all the values of $A$, the computation is repeated keeping the same space for the trial functions, but with a new test space with functions selected only from the three-vortex snapshots. This procedure is repeated similarly for all the possible combinations of trial and test spaces.

\begin{figure}[htbp]
  \centering
  \subfigure[]{
  \includegraphics[width=0.48\textwidth]{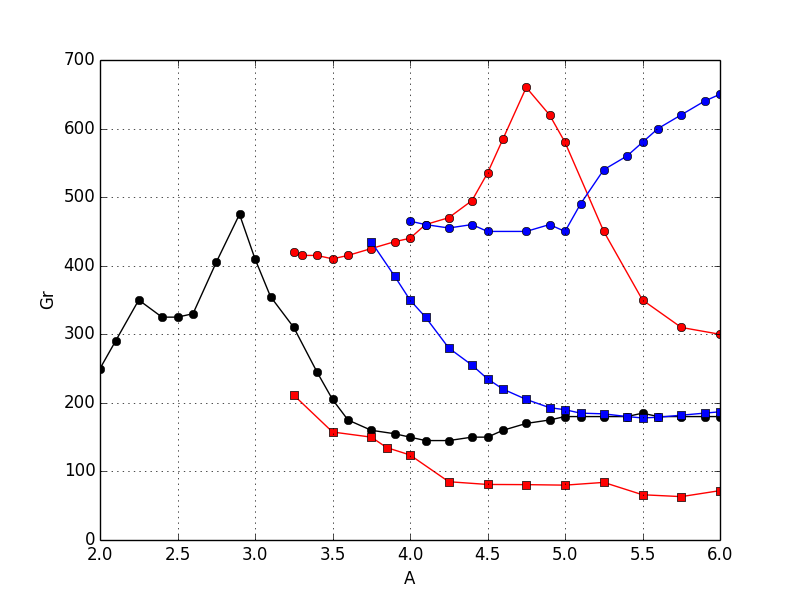}}
  \subfigure[]{
  \includegraphics[width=0.48\textwidth]{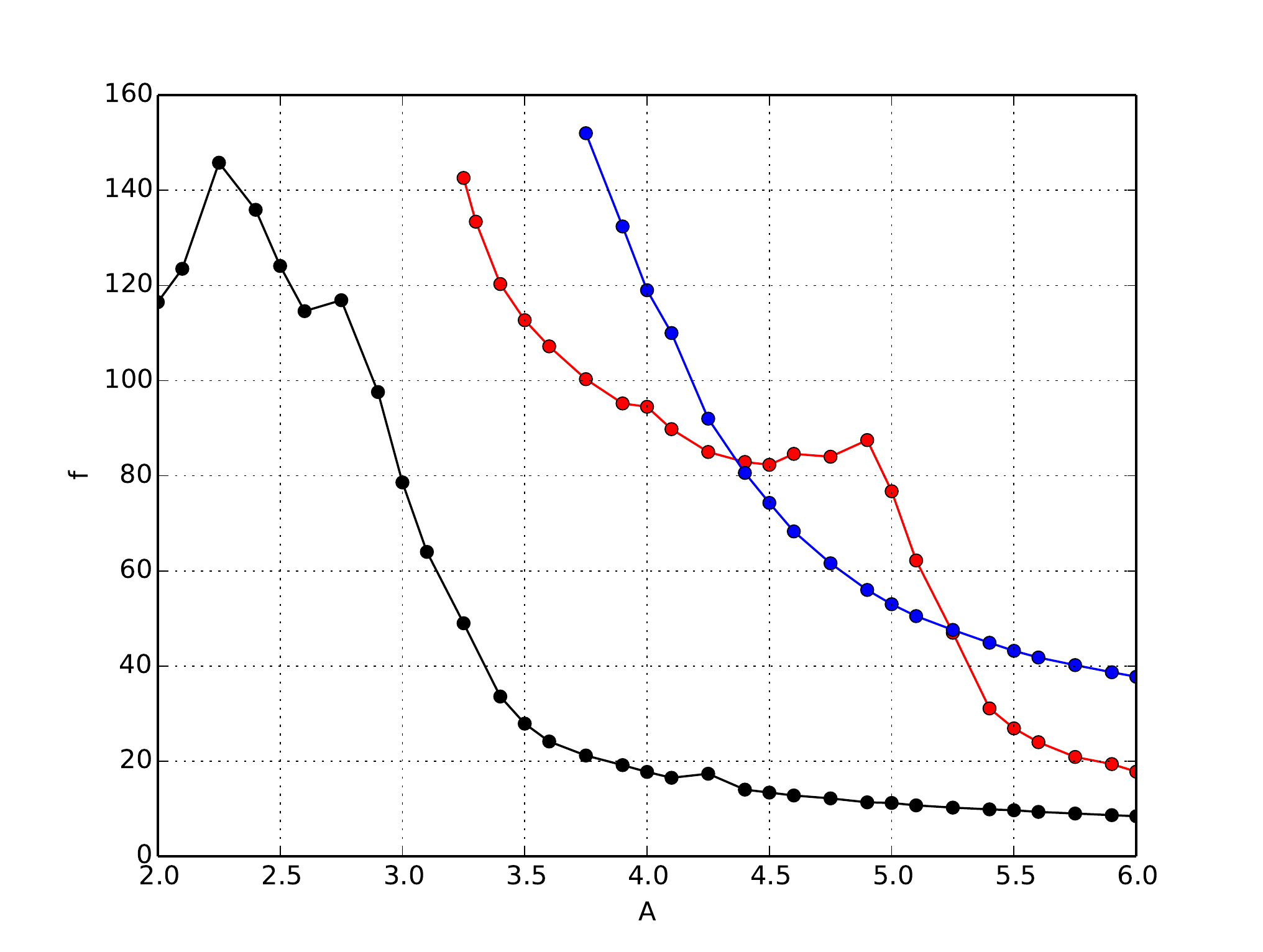}}
  \caption{Left: bifurcation points in the parameters plane for the one-vortex flows (in black), the two-vortex flows (red) and the three-vortex flows (blue). The Hopf bifurcation points are marked with the circles, the steady bifurcation points with the square marks. The Grashof number is expressed in thousands. Right: frequency $f$ at the onset of the periodic solution as a function of the geometric aspect factor $A$. The one-vortex solutions are marked in black, the two-vortex solutions in red and the three-vortex solutions in blue.}
  \label{fig:frequencies}
\end{figure}

The resulting stability regions are shown in Figure~\ref{fig:frequencies} (a), and the oscillation frequency at the onset of the Hopf bifurcation is shown in Figure~\ref{fig:frequencies} (b). We remark that in Figure~\ref{fig:frequencies} (a), the frequency decrease with $A$ follows the benchmark data of~\cite{Gelfgat:Ref11}.

To conclude, we make an estimate of the computational performance of the reduced basis method for the construction of the diagram in Figure~\ref{fig:frequencies} (b).  The reduced order method computed a total of 103 bifurcation points, requiring roughly 43 cpu-hours on a desktop computer. The same computation without using the reduced order method would have taken roughly 5000 cpu-hours. Taking into account the offline phase (the computation of 108 basis functions and the selection of the trial and test spaces), we get the following estimate for the computational time reduction:
\[
  \begin{split}
  \frac{\text{total offline and online computation time}}{\text{equivalent fully offline computation time}}\simeq \\
  \frac{103\cdot 24+20+103\cdot10\cdot\frac{1}{20} \text{ cpu hours}}{103\cdot10\cdot24 \text{ cpu hours}}=10\%,
\end{split}
        \label{eq:savingsestimate}
\]
where the numbers refer to 20 hours for the POD computations, and we consider an average of 10 runs per bifurcation diagram, 103 bifurcations computed, and 5 minutes per run during the online phase and 24 hours per run during the offline phase.

As expected, the computational savings are five times larger than in the one-parameter case, confirming that the proposed RB method becomes more and more convenient as the dimension of the parameter space increases.

\section{Perspectives}
\label{sec:perspectives}
This work highlights several improvements with respect to the state of the art for  stability and bifurcation of parametrized viscous flow problems approached with a reduced order method: in particular the focus has been put on the approximation stability, proper sampling techniques and reduced eigenproblems for stability analysis and bifurcations detection.

On the basis of the experimental results, some remarks concerning the choice of trial and test spaces for the detection of bifurcation points in the reduced order model are made. Without a proper choice of trial and test spaces at the reduced order level, we claim that it is not possible to correctly detect the presence of bifurcation points with the reduced order method.

Numerical tests were performed on a well-known cavity benchmark problem, and where possible the reduced order results are compared with the full order method.


We plan to extend the proposed ROM framework to more advanced 3D problems, in particular with applications to the study of the \emph{Coanda effect} in haemodynamics~\cite{AQpreprint}, and the influence of geometry on symmetry breaking~\cite{PQR_coanda}. More complex studies may be devoted to bifurcations and stability of flows with an elastic thin wall-structure interaction.

\section*{Acknowledgements}
The authors acknowledge Dr. E. Merzari for his help with the \texttt{Nek5000} software and for the useful discussions, and the \texttt{Nek5000} community in general, Dr. F. Ballarin for the insights on approximation stability. G. Pitton has been supported by the pre-doc program at SISSA.
G. Rozza acknowledges the support of NOFYSAS excellence grant program at SISSA and INDAM-GNCS activity group (2015 and 2016 projects), as well as European Union Funding for Research and Innovation – Horizon 2020 Program – in the framework of European Research Council Executive Agency: H2020 ERC CoG 2015 AROMA-CFD project 681447 “Advanced
Reduced Order Methods with Applications in Computational Fluid Dynamics”.
The motivation for developing this work came from Prof. A.T. Patera (MIT) and from Prof. J. Rappaz (EPFL). We acknowledge Prof. F. Brezzi (IUSS, Pavia) for insights and some references.

We gratefully thank Prof. A. Quaini for ongoing collaboration on this topic with the Mathematics Department at University of Houston, USA.

The computing resources have been provided by the Sis14\_COGESTRA cpu time grant allocation at CINECA, Bologna, Italy.

\printbibliography

\end{document}